\documentclass[10pt]{article}
\usepackage{amsmath}
\usepackage{amsfonts}
\usepackage{latexsym}
\usepackage{amssymb}
\usepackage[dvips]{graphics}

\newcommand{\Q}{{\mathbf Q}}
\newcommand{\R}{{\mathbf R}}

\newcommand{\kk}{{\mathbf k}}
\newcommand{\cat}{{\rm {cat }}}
\newcommand{\id}{{\rm {id }}}

\newcommand{\tc}{{\rm\bf {TC}}}

\newtheorem{theorem}{Theorem}

\newtheorem{definition}[theorem]{Definition}

\newenvironment{proof}{{\bf Proof.}}{\par}
\newenvironment{example}{{\bf Example.} }{\par}
\newenvironment{remark}{{\bf Remark.} }{\par}

\title{Topological complexity of motion planning}
\author{Michael Farber\footnote{Partially supported by a grant from the Israel Science Foundation}}
\date{November 18, 2001}

\begin{document}

\sloppy

\maketitle

\begin{abstract}
In this paper we study a notion of topological complexity $\tc(X)$ for the motion planning problem. 
$\tc(X)$ is a number which measures 
discontinuity of the process of motion planning in the configuration space $X$. More precisely, $\tc(X)$ is the minimal number $k$
such that there are $k$ different \lq\lq motion planning rules\rq\rq, 
each defined on an open subset of $X\times X$, so that each rule is
continuous in the source and target configurations. 
We use methods of algebraic topology (the Lusternik - Schnirelman theory) to study the topological complexity $\tc(X)$ .
We give an upper bound for $\tc(X)$ 
(in terms of the dimension of the configuration space $X$) and also a lower bound (in the terms of the structure
of the cohomology
algebra of $X$). 
We explicitly compute the topological complexity of motion planning for a number
of configuration spaces: for spheres, two-dimensional surfaces, for products of spheres.  
In particular, we completely calculate the topological complexity of the problem of motion planning for a robot arm
in the absence of obstacles.\footnote{I am thankful to
D. Halperin, M. Sharir and S. Tabachnikov for a number of very useful conversations.}

\end{abstract}

\section{Definition of topological complexity} 

Let $X$ be the space of all possible configurations of a mechanical system.
In most applications the configuration space $X$ comes equipped with a structure of topological space.
The motion planning problem consists in constructing a program or a devise, which 
takes pairs of configurations $(A,B)\in X\times X$ as an input and produces as an output a continuous path in $X$,
which starts at $A$ and ends at $B$, see \cite{L}, \cite{SS}, \cite{Sh}. Here $A$ is the initial configuration, and 
and $B$ is the final (desired) configuration of the system.

We will assume below that the configuration space $X$ is path-connected, which means that for any pair of points of $X$
there exists a continuous path in $X$ connecting them. Otherwise, the motion planner has first to decide whether the
given points $A$ and $B$ belong to the same path-connected component of $X$.

The motion planning problem can be formalized as follows.
Let $PX$ denote the space of all continuous paths
$\gamma:[0,1]\to X$ in $X$. We will denote by $\pi: PX\to X\times X$ the map associating to any path $\gamma \in PX$ 
the pair of its initial and end points
$\pi(\gamma) =(\gamma(0), \gamma(1))$. Equip the path space $PX$ with the compact-open topology.
Rephrasing the above definition we see that the problem of
motion planning in $X$ consists of finding a function $s: X\times X \to PX$ such that the composition
$\pi\circ s = \id$ is the identity map. In other words, $s$ must be a section of $\pi$.

Does there exist a continuous motion planning in $X$?
Equivalently, we ask whether it is possible to construct a motion planning in the configuration space
$X$ so that the continuous path $s(A,B)$ in $X$, which describes the movement of the system from the initial
configuration $A$
to final configuration $B$, depends continuously on the pair of points $(A, B)$? 
In other words, does there exist a motion planning in $X$ such that the 
section $s: X\times X\to PX$ is continuous? 

Continuity of the motion planning is an important natural requirement.
Absence of continuity will result in instability of the behavior: there will exist arbitrarily close pairs $(A,B)$ and $(A',B')$
of initial - desired configurations such that the corresponding paths $s(A,B)$ and $s(A',B')$ are not close.
\begin{figure}[h]
\begin{center}
\includegraphics[0,0][300,85]{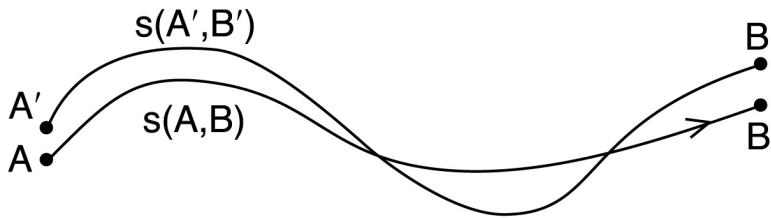}
\end{center}
\caption{Continuity of motion planning: close initial -- final pairs $(A,B)$ and $(A', B')$ produce close movements $s(A,B)$ and $s(A',B')$.}
\end{figure}

Unfortunately, as the following Theorem states, a continuous motion planning exists only in very simple situations.

\begin{theorem}\label{thm1} A continuous motion planning $s:X\times X \to PX$ exists 
if and only if the configuration space $X$
is contractible.
\end{theorem} 
\begin{proof} Suppose that a continuous section $s: X\times X\to PX$ exists. Fix a point $A_0\in X$ and consider
the homotopy 
$$h_t: X\to X, \quad h_t(B)=s(A_0,B)(t),$$ 
where $B\in X$, and $t\in [0,1]$. We have $h_1(B)=B$ and $h_0(B)=A_0$.
Thus $h_t$ gives a contraction of the space $X$ into the point $A_0\in X$. 

Conversely, assume that there is a continuous homotopy $h_t: X\to X$ such that $h_0(A) =A$ and $h_1(A)=A_0$
for any $A\in X$. Given a pair $(A,B)\in X\times X$, we may compose the path $t\mapsto h_t(A)$ with
the inverse of $t\mapsto h_t(B)$, which gives a continuous motion planning in $X$. 

Thus, we get a motion planning in a contractible space $X$ by first
moving $A$ into the base point $A_0$ along the contraction, and then following the inverse of the path, which brings $B$ to
$A_0$. $\diamond$
\end{proof}

\begin{definition} Given a path-connected topological space $X$, we define the {\it topological complexity of the motion planning}
in $X$ as the minimal number ${\rm \bf {TC}}(X)=k$, such that the Cartesian product 
$X\times X$ may be covered by $k$ open subsets
\begin{eqnarray}\label{cover}
X\times X = U_1\cup U_2\cup \dots \cup U_k
\end{eqnarray}
such that for any $i=1, 2, \dots, k$ there exists a continuous motion planning $s_i: U_i\to PX$, $\pi\circ s_i =\id$
over $U_i$. If no such $k$ exists we will set $\tc(X)=\infty$.
\end{definition}

Intuitively, the topological complexity $\tc(X)$ is a measure of discontinuity of any motion planner in $X$.

Given an open cover (\ref{cover}) and sections $s_i$ as above, one may organize a motion planning algorithm as follows. Given a pair
of initial-desired configurations $(A,B)$, we first find the subset $U_i$ with the smallest index $i$ 
such that $(A,B)\in U_i$ and then we give the path $s_i(A,B)$ as an output.
Discontinuity of the output $s_i(A,B)$ as a function of the input $(A,B)$ is obvious: suppose that $(A,B)$ is close to
the boundary of $U_1$ and is close to a pair $(A',B')\in U_2-U_1$; then the output $s_1(A,B)$ compared to 
$s_2(A',B')$ may be completely different, since the sections $s_1|_{U_1\cap U_2}$ and $s_2|_{U_1\cap U_2}$
are in general distinct.
\begin{figure}[h]
\begin{center}
\includegraphics[0,0][260, 210]{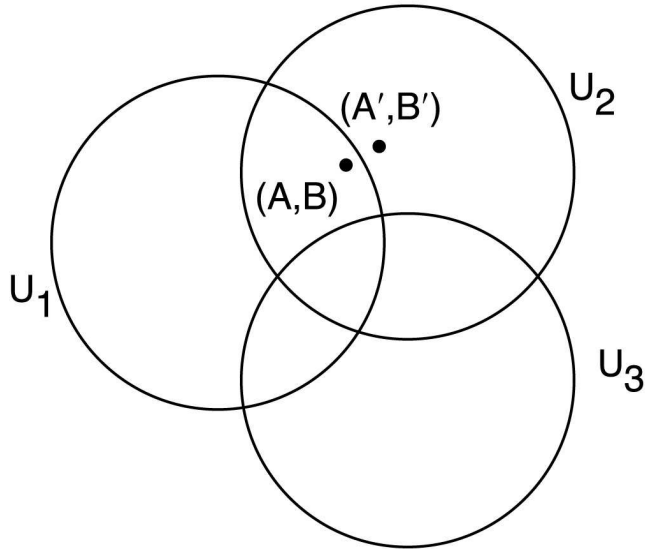}
\end{center}
\caption{Discontinuity of the motion planner corresponding to a covering $\{U_i\}$. }
\end{figure}

According to Theorem \ref{thm1}, we have $\tc(X)=1$ if and only if the space $X$ is contractible.

\begin{example} Suppose that $X$ is a convex subset of an Euclidean space $\R^n$. Given a pair of 
initial - desired configurations
$(A,B)$, we may move with constant velocity 
along the straight line segment connecting $A$ and $B$. This clearly produces a continuous algorithm for the
motion planning problem in $X$. This is consistent with Theorem \ref{thm1}: we have $\tc(X)=1$ since $X$ is contractible.
\end{example}

\begin{example}
Consider the case when $X=S^1$ is a circle. Since $S^1$ is not contractible, we know that
$\tc(S^1)>1$. Let us show that ${\tc}(S^1)=2$. 
Define $U_1\subset S^1\times S^1$ as $U_1=\{(A,B); A\not= -B\}$. A continuous motion planning over
$U_1$ is given by the map $s_1: U_1\to PS^1$ which moves $A$ towards $B$ with constant velocity along the 
unique shortest arc connecting $A$ to $B$. This map $s_1$ cannot be extended to a continuous map on the 
pairs of antipodal points $A=-B$. 
Now define $U_2=\{(A,B); A\not=B\}$. Fix an orientation of the circle $S^1$. A continuous motion planning over
$U_2$ is given by the map $s_2: U_2\to PS^1$ which moves $A$ towards $B$ with constant velocity in the positive 
direction along the circle. Again, $s_2$ cannot be extended to a continuous map on the whole $S^1\times S^1$. 
\end{example}

\begin{remark} Our definition of the topological complexity $\tc(X)$ is motivated by the notion of a genus of a fiber space, introduced
by A.S. Schwarz \cite{Sz}. In fact $\tc(X)$ is the Schwarz genus of the path space fibration $PX\to X\times X$.

The theory of Schwarz genus was used by S. Smale \cite{Sm1} and V. Vassiliev \cite{V1}, \cite{V2} 
to define topological complexity
of algorithms of finding roots of polynomial equations. 
\end{remark} 

\section{Homotopy invariance}

The following property of homotopy invariance allows often to simplify the configuration space $X$ without changing the 
topological complexity $\tc(X)$.

\begin{theorem}\label{thm2}  $\tc(X)$ 
depends only on the homotopy type of $X$.
\end{theorem}
\begin{proof} Suppose that $X$ {\it dominates} $Y$, i.e. there exist continuous maps $f: X\to Y$ and $g:Y\to X$
such that $f\circ g\simeq \id_Y$. Let us show that then $\tc(Y)\leq \tc(X)$. 
Assume that $U\subset X\times X$ is an open subset such that there exists a continuous motion planning 
$s: U\to PX$ over $U$. Define $V=(g\times g)^{-1}(U)\subset Y\times Y$. We will construct a continuous motion
planning $\sigma: V\to PY$ over $V$ explicitly. 
Fix a homotopy $h_t: Y\to Y$ with $h_0=\id_Y$ and $h_1=f\circ g$; here $t\in [0,1]$. 
For $(A,B)\in V$ and $\tau\in [0,1]$ set
\begin{eqnarray*}
\sigma(A,B)(\tau) = 
\left\{
\begin{array}{ll}
h_{3\tau}(A), &\mbox{for}\quad 0\leq \tau\leq 1/3, \\ \\
f(s(gA,gB)(3\tau-1)),&\mbox{for}\quad 1/3\leq \tau\leq 2/3, \\ \\
h_{3(1-\tau)},&\mbox{for}\quad 2/3\leq \tau\leq 1.
\end{array}
\right.
\end{eqnarray*}
Thus we obtain that for $k=\tc(X)$ any open cover $U_1\cup\dots \cup U_k=X\times X$ with a continuous motion planning 
over each $U_i$ defines an open cover $V_1\cup\dots \cup V_k$ of $Y\times Y$
with the similar properties. This proves that $\tc(Y)\leq \tc(X)$, and obviously implies the statement 
of the Theorem. $\diamond$
\end{proof}

\section{An upper bound for $\tc(X)$}

\begin{theorem}\label{cor1} For any path-connected paracompact space $X$, we have 
\begin{eqnarray}
\tc(X)\, \leq \, 2\cdot \dim X +1.\label{ineq1}
\end{eqnarray}
In particular, if $X$ is a connected polyhedral subset of $\R^n$ then the topological complexity
$\tc(X)$ can be estimated from above as follows
\begin{eqnarray}
\tc(X) \, \leq \, 2n-1.\label{ineq2}
\end{eqnarray}
\end{theorem}

We postpone the proof.

We will use a relation between $\tc(X)$ and the Lusternik - Schni\-relman category $\cat(X)$.
Recall that $\cat(X)$ is defined as the smallest integer $k$ such that $X$ may be covered by $k$ open subsets
$V_1\cup \dots \cup V_k=X$ with each inclusion $V_i\to X$ null-homotopic.

\begin{theorem}\label{thm3} If $X$ is path-connected and paracompact then
\begin{eqnarray}
\cat(X)\, \leq\,  \tc(X)\, \leq\,  2\cdot \cat(X)-1.\label{ineq3}
\end{eqnarray}
\end{theorem}
\begin{proof} Let $U\subset X\times X$ be an open subset such that there exists a continuous motion planning
$s: U\to PX$ over $U$. Let $A_0\in X$ be a fixed point. Denote by $V\subset X$ the set of all points $B\in X$
such that $(A_0,B)$ belongs to $U$. Then clearly the set $V$ is open and it is contractible in $X\times X$.

If $\tc(X)=k$ and $U_1\cup\dots\cup U_k$ is a covering of $X\times X$ with a continuous motion planning over 
each $U_i$, then the sets $V_i$, where $A_0\times V_i = U_i\cap (A_0\times X)$ form a categorical open cover of $X$.
This shows that $\tc(X)\geq \cat(X)$. 

The second inequality follows from the obvious inequality
$$\tc(X) \leq \cat(X\times X)$$ 
combined with $\cat(X\times X)\leq 2\cdot \cat(X) -1$, see Proposition 2.3 of \cite{J}. $\diamond$
\end{proof} 

\vskip 0.5cm
\noindent
{\bf Proof of Theorem \ref{cor1}}. It is well-known that $\cat(X)\leq \dim (X) +1$. Together with the right
inequality in (\ref{ineq3}) this gives (\ref{ineq1}). 

If $X\subset \R^n$ is a connected polyhedral subset then $X$ has homotopy type of an $(n-1)$-dimensional polyhedron $Y$.
Using homotopy invariance (Theorem \ref{thm2}) we find $\tc(X)=\tc(Y)\leq 2(n-1)+1=2n-1$.
$\diamond$

\section{A lower bound for $\tc(X)$}

Let $\kk$ be a field. The cohomology $H^\ast(X;\kk)$ is a graded $\kk$-algebra with the multiplication
\begin{eqnarray}
\cup: H^\ast(X;\kk)\otimes H^\ast(X;\kk)\to H^\ast(X;\kk)\label{prod}
\end{eqnarray}
given by the cup-product, see \cite{DNF}, \cite{Sp}. 
The tensor product $H^\ast(X;\kk)\otimes H^\ast(X;\kk)$ is also a graded $\kk$-algebra
with the multiplication 
\begin{eqnarray}\label{signs}
(u_1\otimes v_1)\cdot (u_2\otimes v_2) = (-1)^{|v_1|\cdot |u_2|}\, u_1u_2\otimes v_1v_2.
\end{eqnarray}
Here $|v_1|$ and $|u_2|$ denote the degrees of cohomology classes $v_1$ and $u_2$ correspondingly.
The cup-product (\ref{prod}) is an algebra homomorphism.

\begin{definition} The kernel of homomorphism (\ref{prod}) will be called {\it the ideal of zero-divisors} of $H^\ast(X;\kk)$.
The {\it zero-divisors-cup-length} of $H^\ast(X;\kk)$ is the length of the longest nontrivial product in the ideals of zero-divisors
of $H^\ast(X;\kk)$.
\end{definition}
\begin{example} Let $X=S^n$. Let $u\in H^n(S^n;\kk)$
be the fundamental class, and let $1\in H^0(S^n;\kk)$ be the unit.
Then $a=1\otimes u-u\otimes 1\in H^\ast(S^n;\kk)\otimes H^\ast(S^n;\kk)$ is a zero-divisor, since
applying homomorphism (\ref{prod}) to $a$ we obtain 
$ 1\cdot u-u\cdot 1=0.$ 
Another zero-divisor is $b=u\otimes u$, since $b^2=0$. 
Computing $a^2=a\cdot a$ by means of rule (\ref{signs}) we find
$$a^2=((-1)^{n-1}-1)\cdot u\otimes u.$$ 
Hence $a^2=-2 b$ for $n$ even and $a^2=0$ for $n$ odd;
the product $ab$ vanishes for any $n$.
We conclude that {\it the zero-divisors-cup-length of $H^\ast(S^n;\Q)$
equals 1 for $n$ odd and 2 for $n$ even.} \end{example}

\begin{theorem}\label{thm4} The topological complexity of motion planning $\tc(X)$ 
is greater than the zero-divisors-cup-length of 
$H^\ast(X;\kk)$.
\end{theorem}

To illustrate this Theorem, consider the special case $X=S^n$.
Using the computation of the zero-divisors-cup-length for $S^n$
(see the example above) and applying Theorem \ref{thm4}
we find that $\tc(S^n)>1$ for $n$ odd and $\tc(S^n)>2$ for $n$ even.
This means that any motion planner on the sphere $S^n$ must have at least two open sets $U_i$; 
moreover, any motion planner on the sphere $S^n$ must have at least 
three open sets $U_i$ if $n$ is even.

\noindent
\begin{proof} Consider the following commutative diagram
\begin{eqnarray*}
\begin{array}{ccc}
X & \stackrel \alpha \to & PX\\ 
& \stackrel\searrow \Delta & \downarrow \pi\\ 
& & X\times X
\end{array}
\end{eqnarray*}
Here $\alpha$ associates to any point $x\in X$ the constant path $[0,1]\to X$ at this point. 
$\Delta: X\to X\times X$
is the diagonal map $\Delta(x)=(x,x)$. Note that $\alpha$ is a homotopy equivalence. 
The composition
\begin{eqnarray}
H^\ast(X;\kk)\otimes H^\ast(X;\kk) \simeq H^\ast(X\times X;\kk)\stackrel {\pi^\ast}\to H^\ast(PX;\kk)
\underset{\simeq}{\stackrel {\alpha^\ast}\to} H^\ast(X;\kk)\label{isom}
\end{eqnarray} 
coincides with
the cup-product homomorphism (\ref{prod}). Here the homomorphism on the left is the K\"unneth isomorphism.

As we mentioned above, the topological complexity of motion planning $\tc(X)$ is the Schwarz genus (cf. \cite{Sz})
of the fibration $\pi: PX\to X\times X$. 
The statement of Theorem \ref{thm4}
follows from our remarks above concerning homomorphism (\ref{isom}) and from the 
cohomological lower bound for the Schwarz genus, see Theorem 4 of \cite{Sz}.
$\diamond$
\end{proof}

\section{Motion planning on spheres}

\begin{theorem}\label{thm5} The topological complexity of motion planning on the $n$-dimen\-sional sphere $S^n$
is given by
\begin{eqnarray*}
\tc(S^n) \, =\, \left\{
\begin{array}{ll}
2, &\mbox{for $n$ odd},\\ \\
3, &\mbox{for $n$ even.}
\end{array}
\right.
\end{eqnarray*}
\end{theorem}
\begin{proof} First we will show that $\tc(S^n)\leq 2$ for $n$ odd. 
Let $U_1\subset S^n\times S^n$ be the set of all pairs $(A,B)$ where $A\not=-B$.
Then there is a unique shortest arc of $S^n$ connecting $A$ and $B$ and we will construct a continuous motion planning
$s_1: U_1\to PS^n$ by setting $s_1(A,B)\in PS^n$ to be this shortest arc passed with a constant velocity.
The second open set will be defined as $U_2=\{(A,B); A\not=B\}\subset S^n\times S^n$. A continuous motion 
planning over $U_2$ will be constructed in two steps. On the first step we will move the initial point $A$ to the antipodal
point $-B$ along the shortest arc as above. On the second step we will move the antipodal point $-B$ to $B$. For this purpose
fix a continuous tangent vector field $v$ on $S^n$, which is nonzero
at every point; here we will
use the assumption that the dimension $n$ is odd. We may move $-B$ to $B$ along the spherical arc
$$-\cos \pi t\cdot B + \sin \pi t\cdot v(B), \quad t\in [0,1].$$
This proves that $\tc(S^n)\leq 2$ for $n$ odd; hence by Theorem \ref{thm1} $\tc(S^n)=2$ for $n$ odd.

Assume now that $n$ is even. Let us show that then $\tc(S^n)\leq 3$.  
We will define a continuous motion planning over a set $U_1\subset S^n\times S^n$ as above.
For $n$ even we may construct a continuous tangent vector field $v$ on $S^n$, which vanishes at a single point
$B_0\in S^n$ and is nonzero for any $B\in S^n$, $B\not= B_0$. We will define the second set 
$U_2\subset S^n\times S^n$ as $\{(A,B); A\not= B \, \& \,  B\not= B_0\}$. We may define $s_2: U_2\to PS^n$ as above.
Now, $U_1\cup U_2$ covers everything except the pair of points $(-B_0, B_0)$. Chose a point $C\in S^n$, distinct from $B_0, -B_0$
and set $U_3= S^n-C$. Note that $U_3$ is diffeomorphic to $\R^n$ and so there exists a continuous motion planning over $U_3$.
This proves that $\tc(S^n)\leq 3$. On the other hand, using Theorem \ref{thm4} and the preceeding Example, we find $\tc(S^n)\geq 3$
for $n$ even. This completes the proof.
$\diamond$
\end{proof}
\section{More examples}

\begin{theorem} Let $X=\Sigma_g$ be a compact orientable two-dimensional surface of genus $g$. Then 
\begin{eqnarray*}
\tc(X) \, =\, \left\{
\begin{array}{lll}
3, & \mbox{if}& g\le 1,\\
5, & \mbox{if}& g>1.
\end{array}
\right.
\end{eqnarray*}
\end{theorem} Consider first the case $g\geq 2 $. Then we may find cohomology classes $u_1, v_1, u_2, v_2\in H^1(X;\Q)$
forming a symplectic system, i.e. $u_i^2=0$, $v_i^2=0$, and 
$u_1v_1=u_2v_2=A\not=0$, where
$A\in H^2(\Sigma_g;\Q)$ is the fundamental class; besides,
$v_iu_j=v_iv_j=u_iu_j=0$ for $i\not= j$.
Then in the algebra $H^\ast(X;\Q)\otimes H^\ast(X;\Q)$ holds
$$\prod_{i=1}^2 (1\otimes u_i-u_i\otimes 1)(1\otimes v_i-v_i\otimes 1)= 2A\otimes A\, \not= 0$$
and hence we obtain, using Theorem \ref{thm4}, that $\tc(X)\ge 5$. The opposite inequality follows from Theorem \ref {cor1}.

The case $g=0$ follows from Theorem \ref{thm5} since then $X=S^2$. The case $g=1$, which corresponds to the two-dimensional torus $T^2$,
will be considered later in Theorem \ref{robot}. $\diamond$

\begin{theorem} Let $X={\mathbf {CP}}^n$ be the $n$-dimensional complex projective space. Then $\tc(X)\geq 2n+1$.
\end{theorem}
\begin{proof} If $u\in H^2(X;\Q)$ is a generator, then 
\begin{eqnarray*}
(1\otimes u-u\otimes 1)^{2n}\, =\, (-1)^n 
\left(
\begin{array}{c}
 2n \\
n
\end{array}
\right)
 u^n\otimes u^n\not= 0.
\end{eqnarray*}
Hence Theorem \ref{thm4} gives $\tc(X)\ge 2n+1$. $\diamond$
\end{proof}

\section{Product inequality}

\begin{theorem}\label{prod1} For any path--connected metric spaces
$X$ and $Y$,
\begin{eqnarray}
\tc(X\times Y)\leq \tc(X) +\tc(Y) -1.
\end{eqnarray}
\end{theorem}
\begin{proof} Denote $\tc(X)=n$, $\tc(Y)=m$. Let $U_1,\dots, U_n$ be on open cover of $X\times X$ with a continuous motion 
planning $s_i: U_i\to PX$ for $i=1, \dots, n$.  Let $f_i: X\times X\to \R$, where $i=1, \dots, n$, be a partition of unity
subordinate to the cover $\{U_i\}$.
Similarly, 
let $V_1,\dots, V_m$ be on open cover of $Y\times Y$ with a continuous motion 
planning $\sigma_j: V_j\to PY$ for $j=1, \dots, m$, and let $g_j: Y\times Y\to \R$, where $j=1, \dots, m$ be a partition of unity subordinate to the cover $\{V_j\}$. 

For any pair of nonempty subsets $S\subset \{1, \dots, n\}$ and $T\subset \{1, \dots, m\}$, let 
$$W(S,T)\subset (X\times Y)\times (X\times Y)$$ 
denote the set of all
4-tuples $(A,B, C,D)\in (X\times Y)\times (X\times Y)$, such that for any $(i,j)\in S\times T$ and for any $(i',j')\notin S\times T$ holds
$$f_i(A,C) \cdot g_j(B,D) > f_{i'}(A,C) \cdot g_{j'}(B,D).$$
One easily checks that:

{\it (a) each set $W(S,T)\subset X\times X$ is open;

(b) $W(S,T)$ and $W(S', T')$ are 
disjoint if neither $S\times T\subset S'\times T'$, nor $S'\times T'\subset S\times T$;

(c) if $(i,j)\in S\times T$, then $W(S,T)$ is contained in $U_i\times V_j$; therefore there exists a continuous motion planning 
over each $W(S,T)$ (it can be described explicitly in terms of $s_i$ and $\sigma_j$);

(d) the sets $W(S,T)$ (with all possible nonempty $S$ and $T$) form a cover of $(X\times Y)\times (X\times Y)$.}

Let us prove (d). Suppose that $(A,B,C,D)\in (X\times Y)\times (X\times Y)$. Let $S$ be the set of all indices 
$i\in \{1, \dots, n\}$,
such that $f_i(A,C)$ equals the maximum of $f_k(A,C)$, where $k=1, 2, \dots, n$. 
Similarly, let $T$ be the set of all $j\in \{1, \dots, m\}$, such that $g_j(B,D)$ equals the maximum of $g_\ell(B,C)$, where
$\ell=1, \dots, m$. Then clearly $(A,B,C,D)$ belongs to $W(S,T)$.

Let $W_k \subset (X\times Y)\times (X\times Y)$ 
denote the union of all sets $W(S,T)$, where $|S|+|T|=k$. Here $k=2, 3, \dots, n+m.$
The sets $W_2, \dots, W_{n+m}$ form an open cover of $(X\times Y)\times (X\times Y)$. 
If $|S|+|T|=|S'|+|T|=k,$ then the corresponding sets $W(S,T)$ and $W(S',T')$ either coincide (if $S=S'$ and $T=T'$), or are disjoint. 
Hence we see (using (c)) that there exists a continuous motion planning over each open set $W_k$. 

This completes the proof.
 $\diamond$
\end{proof}

\begin{remark} The above proof represents a modification of the arguments of the proof of the product
inequality for the Lusternik - Schnirelman category, see page 333 of \cite{J}.
\end{remark}

\section{Motion planning for a robot arm}

Consider a robot arm consisting of $n$ bars $L_1, \dots, L_n$, such that $L_i$ and $L_{i+1}$ are connected by
flexible joins. We assume that the initial point of $L_1$ is fixed. In the planar case, a configuration of the arm is determined
by $n$ angles $\alpha_1, \dots, \alpha_n$, where $\alpha_i$ is the angle between $L_i$ and the $x$-axis. 
Thus, in the planar case,
the configuration space of the robot arm (when no obstacles are present) is the $n$-dimensional torus 
$$T^n=S^1\times S^1\times \dots \times S^1.$$
Similarly, the configuration space of a robot arm in the 3-dimensional space $\R^3$ is the Cartesian product of $n$ copies of the two-dimensional sphere $S^2$. 

\begin{figure}[h]
\begin{center}
\includegraphics[0,0][216,128]{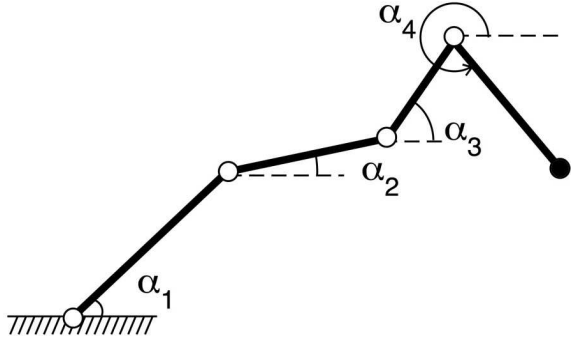}
\end{center}
\caption{Planar robot arm}
\end{figure}

\begin{theorem}\label{arm} The topological complexity of motion planning problem 
of a plane $n$-bar robot arm 
equals $n+1$.  The topological complexity of motion planning problem 
of a spacial $n$-bar robot arm equals $2n+1$. 
\end{theorem}

\begin{remark} It is not difficult to explicitly construct motion planners for the planar and spacial robot arms, which have
the minimal possible topological complexity. 
Such algorithms could be based on the ideas used in the proof of the product
inequality (Theorem \ref{prod1}).
\end{remark}

Theorem \ref{arm} automatically follows from the next statement:

\begin{theorem}\label{robot} Let $X=S^m\times S^m\times \dots \times S^m$ 
be a Cartesian product of $n$ copies of the $m$-dimensional sphere $S^m$. Then 
\begin{eqnarray}\label{stam22}
\tc(X) \, =\, 
\left\{
\begin{array}{ll}
n+1, & \mbox{if $m$ is odd},\\
2n+1, & \mbox{if $m$ is even}.
\end{array}
\right.
\end{eqnarray}
\end{theorem}
\begin{proof} Using the product inequality (Theorem \ref{prod1}) and the calculation for spheres (Theorem \ref{thm5})
we find that $\tc(X)$ is less or equal than the RHS of (\ref{stam22}).
To establish the inverse inequality we will use Theorem \ref{thm4}.
Let $a_i\in H^m(X;\Q)$ denote the cohomology class which is the pull-back of the fundamental class of $S^m$ under the
projection $X\to S^m$ onto the $i$-th factor; here $i=1, 2, \dots, n$. We see that 
$$\prod_{i=1}^n (1\otimes a_i-a_i\otimes 1)\not=\, 0\in H^\ast(X\times X;\Q).$$
This shows that the zero-divisors-cup-length of $X$ is at least $n$. If $m$ is even then 
$$\prod_{i=1}^n (1\otimes a_i-a_i\otimes 1)^2\not=\, 0\in H^\ast(X\times X;\Q).$$
Hence for $m$ even, the zero-divisors-cup-length of $X$ is at least $2n$. 
Application of Theorem \ref{thm4} completes the proof. $\diamond$
\end{proof}

\bibliographystyle{amsalpha}

\vskip 2cm 

Address: 

Michael Farber,

School of Mathematical Sciences, 

Tel Aviv University, Ramat Aviv 69978, Israel

farber@math.tau.ac.il

\end{document}